\documentclass[12pt]{amsart}
\usepackage[utf8]{inputenc}

\usepackage{amsmath,amsthm,amssymb,mathscinet,bbm}
\usepackage{parskip}
\usepackage{geometry}
\usepackage{xcolor} 
\usepackage{xfrac, nicefrac}

\usepackage{mathtools} 
\usepackage[hidelinks]{hyperref}
\author{Paul Pollack} 
\author{Akash Singha Roy}
\address{Department of Mathematics \\ University of Georgia \\ Athens, GA 30602}
\email{pollack@uga.edu}
\email{akash01s.roy@gmail.com}


\renewcommand\phi\varphi
\usepackage{accents}

\renewcommand{\pod}[1]{\allowbreak\mathchoice
  {\if@display \mkern 18mu\else \mkern 8mu\fi (#1)}
  {\if@display \mkern 18mu\else \mkern 8mu\fi (#1)}
  {\mkern4mu(#1)}
  {\mkern4mu(#1)}
}
\usepackage{graphicx}
\DeclareMathAlphabet{\curly}{U}{rsfs}{m}{n}

\newcommand{\1}{\mathbbm{1}}

\newcommand\NatNos{\mathbb N}

\newtheorem{thm}{Theorem}[section]

\newtheorem{prop}[thm]{Proposition}
\newtheorem{lem}[thm]{Lemma}

\theoremstyle{remark}
\newtheorem*{rmk}{Remark}

\newcommand\err{\mathcal E}
\newcommand\reals{\mathbb R}
\newcommand\complex{\mathbb C}
\newcommand\Udisk{\mathbb U}
\newcommand\sumxyzf{\sum\nolimits_{n \le x}^* f(n)}
\newcommand\expfpp{\exp\left(\sum_{p \le y} \frac{|f(p)|}p \right)}

\newcommand\summ{\sum_{\substack{m \le x\\P^+(m) \le y}}}

\newcommand\sumj{\sum_{j \ge 1}}
\newcommand\sumPTwoj{\sum_{\substack{P_2, \dots , P_j \in (y, P_1)\\P_2, \dots, P_j \text{ distinct}\\P_2 \cdots P_j \le x/mP_1}}}
\newcommand\sumPTwojMO{\sum_{\substack{P_2, \dots, P_{j-1} \in (y, P_1)\\P_2, \dots, P_{j-1} \text{ distinct}\\P_2 \cdots P_{j-1} \le x/myP_1}}}

\newcommand\bbm{\mathbbm 1}

\newcommand\chibara{\overline\chi(a)}

\newcommand\chiphin{\chi(\phi(n))}
\newcommand\sumnontrivchar{\sum_{\chi \ne \chi_0 \bmod q}}
\newcommand\phiqrec{\frac1{\phi(q)}}
\newcommand\condchi{\mathfrak f(\chi)}
\newcommand\cond{\mathfrak f}
\newcommand\quadxtoinfty{\quad\text{as $x\to\infty$}}

\numberwithin{equation}{section}
\begin{document}
\title[Mean Values of Multiplicative Functions and residue-class distribution]{Mean values of multiplicative functions and applications to residue-class distribution}

\begin{abstract} We provide a uniform bound on the partial sums of multiplicative functions under very general hypotheses. As an application, we give a nearly optimal estimate for the count of $n \le x$ for which the Alladi-Erd\H{o}s function $A(n) = \sum_{p^k \parallel n} k p$ takes values in a given residue class modulo $q$, where $q$ varies uniformly up to a fixed power of $\log x$. We establish a similar result for the equidistribution of the Euler totient function $\phi(n)$ among the coprime residues to the ``correct" moduli $q$ that vary uniformly in a similar range, and also quantify the failure of  equidistribution of the values of $\phi(n)$ among the coprime residue classes to the ``incorrect" moduli. 
\end{abstract}

\maketitle

\section{Introduction}
The subject of mean values of multiplicative functions has received a lot of attention in the field of analytic number theory. The aim is to provide estimates for the partial sums of multiplicative functions that are as precise and as general as possible. In general settings, a satisfactory bound is often provided by the classical result of Hal\'asz \cite{halasz68}, made quantitatively precise in work of Montgomery and Tenenbaum (see Corollary 4.12 on p. 494 of \cite{tenenbaum15} for the precise statement); several useful variants of this result may also be found in \cite{HT88}, \cite{HT91} and  \cite{tenenbaum15}. In the direction of precision, one of the most powerful estimates is provided by the method of Landau--Selberg--Delange, a comprehensive account of which may be found in \cite[Chapter II.5]{tenenbaum15}. However, while the estimate of Hal\'asz is in certain occasions not precise enough for applications, the error terms arising from the Landau--Selberg--Delange method are sometimes too large to give useful estimates in practical settings. One of our primary objectives in this manuscript is to bridge this gap between precision and generality, so as to obtain a bound that captures certain power savings missed by Hal\'asz's Theorem in applications, whilst having error terms that are small enough to permit applicability in a wide variety of settings. In order to demonstrate the flexibility of both our result and our methods, we will obtain precise estimates on the value distributions of the Alladi-Erd\H{o}s function $A(n) \coloneqq \sum_{p^k \parallel n} k p$ and the Euler totient function $\phi(n) \coloneqq \#\{1 \le d \le n: ~ \gcd(d, n)=1 \}$ in coprime residue classes to moduli varying uniformly in a wide range.     

In what follows, $\Udisk$ denotes the unit disk in the complex plane, namely, the set $\{s \in \complex: |s| \le 1\}$. For given real number $z$, we say that a positive integer $n$ is $z$\textsf{-smooth} if it has no prime factor exceeding $z$. The number of $z$-smooth numbers up to $x$ is denoted by $\Psi(x, z)$, a quantity which has been studied extensively in the literature.   
\begin{thm}\label{thm:mainbound}
Let $f \colon \NatNos \rightarrow \Udisk$ be a multiplicative function and $x, y, z, M$ be positive real numbers such that $M \ge 1$, $4<y \le z^{1/2}$ and $z<x$. Assume that there exists a complex number $\varrho \in \mathbb U$ such that for all $Y \ge y$,
\begin{equation}\label{eq:mainhypo}
\sum_{y<p \le Y} f(p) = \varrho (\pi(Y) - \pi(y)) + O(M Y \err(y)), \end{equation}
for some decreasing function $\err: \reals^+ \rightarrow \reals^+$ satisfying $\lim_{X \rightarrow \infty} \err(X) = 0$. Then
$$\sum_{n \le x} f(n) \ll \frac{|\varrho|x}{\log z} \left(\frac{\log x}{\log y}\right)^{|\varrho|} \expfpp + \Psi(x, z) + Mx (\log x)^2 \left(\err(y)+\frac1y\right),$$
where the implied constant depends at most on the implied constant in \eqref{eq:mainhypo}. 
\end{thm}
In applications, the terms involving $\Psi(x, z)$, $\err(y)$ and $1/y$ often become very small, so the usefulness of the bound is dictated by the size of the first expression. 
The fact that this expression appears with a factor of $|\varrho|$ permits a wide applicability of Theorem \ref{thm:mainbound}. 

Let $A(n) \coloneqq \sum_{p^k\parallel n} kp$ denote the sum of the prime divisors of $n$ counted with multiplicity, known as the \textsf{Alladi-Erd\H{o}s function}. 
A widely applicable criterion of Delange \cite{delange69} gives necessary and sufficient conditions for an integer-valued additive function to be equidistributed among the residue classes to a fixed modulus. This criterion shows that $A(n)$ is equidistributed modulo any fixed positive integer $q$. More precise results were established by Alladi and Erdos \cite{alladi77} for $q=2$ and by Goldfeld \cite{goldfeld17} for general fixed $q$. In recent work \cite{PSR23}, we proved the equidistribution of $A(n)$ mod $q$, uniformly for $q \le (\log x)^{K}$; in other words, we established that  
$$\#\{n \le x: ~ A(n) \equiv a \pmod q\} \sim \frac xq \quadxtoinfty,$$
uniformly in $q \le (\log x)^K$ and in residue classes $a$ mod $q$ (see \cite{PSR22A} for our earlier, weaker result using different methods). Goldfeld noted in \cite{goldfeld17} that for fixed $q$, the deviation from equidistribution is always $O_q(x/(\log x)^{1/2})$. We establish a uniform version of Goldfeld's result. 
\begin{thm}\label{thm:A(n)Effective}
Fix $K \ge 1$ and $\epsilon \in (0, 1)$. We have
$$\#\{n \le x: A(n) \equiv a \pmod q\} = \frac xq + O_{K, \epsilon}\left(\frac x{q(\log x)^{1/2-\epsilon}}\right),$$
uniformly in $q \le (\log x)^K$ and residue classes $a$ mod $q$. 
\end{thm} 
It follows from Goldfeld's methods in \cite{goldfeld17} that for $q=6$, there is a secondary term of order $x/(\log x)^{1/2}$, which means that the exponent $1/2$ in Theorem \ref{thm:A(n)Effective} is essentially best possible. 
Our arguments for Theorem \ref{thm:mainbound} are flexible and can also be modified to yield uniform bounds on the partial sums of multiplicative functions, when the inputs $n$ are restricted to those with sufficiently many large prime factors. Such estimates turn out to be very useful in studying the asymptotic distribution of some classical multiplicative functions, such as Euler's totient function $\phi(n)$, in coprime residue classes to moduli varying uniformly in a wide range. 

Following Narkiewicz  \cite{narkiewicz66}, we say that an integer-valued arithmetic function $g$ is \textsf{weakly uniformly distributed} (or \textsf{weakly equidistributed}) modulo $q$ if $\gcd(g(n),q)=1$ for infinitely many $n$ and, for every \emph{coprime} residue class $a\bmod{q}$,
\begin{equation}\label{eq:WUDrelation} \#\{n\le x: g(n)\equiv a\pmod{q}\} \sim \frac{1}{\phi(q)} {\#\{n\le x: \gcd(g(n),q)=1\}} , \quad\text{as $x\to\infty$}. \end{equation}
Narkiewicz shows in \cite{narkiewicz66} that $\phi(n)$ is weakly equidistributed precisely to those moduli that are coprime to $6$. Dence and Pomerance \cite{DP98} also give asymptotic formulae measuring the failure of weak equidistribution of $\phi(n)$ mod $3$. 

It seems of some of interest to establish weak equidistribution results without the restriction to a fixed modulus. Here we say that 
an integer-valued arithmetic function $g(n)$ is \textsf{weakly equidistributed} mod $q$, \textsf{uniformly for $q \le (\log x)^{K}$,} if:
\begin{enumerate}
\item[(i)] For every such $q$, $g(n)$ is coprime to $q$ for infinitely many $n$, and

\item[(ii)] The relation \eqref{eq:WUDrelation} holds uniformly in moduli $q \le (\log x)^K$ and in coprime residue classes $a$ mod $q$. Explicitly, this means that for any $\epsilon>0$, there exists $X(\epsilon)>0$ such that the ratio of the left hand side of \eqref{eq:WUDrelation} to the right hand side lies in $(1-\epsilon, 1+\epsilon)$ for all $x>X(\epsilon)$, $q \le (\log x)^K$ and coprime residues $a$ mod $q$.  
\end{enumerate}
The weak equidistribution of certain classes of arithmetic functions to uniformly varying moduli was first investigated by the authors in \cite{LPR21, PSR22}. In \cite{PSR23}, we prove a general theorem guaranteeing the weak equidistribution of $\phi(n)$ uniformly to moduli $q \le (\log x)^K$ coprime to $6$; the key idea is to exploit a mixing phenomenon within the group of units mod $q$ (see also related ideas of De Koninck and K\'atai in \cite{DKK06} and \cite{Kat07}). While it is not worthwhile to study the distribution of $\phi(n)$ among the coprime residue classes to even moduli, the corresponding distribution to odd moduli divisible by $3$ is still an interesting question left unaddressed by the aforementioned results. Furthermore, the error terms that arise from the arguments in \cite{PSR23} are quite weak. Our final two theorems address these defects. We first establish an effective estimate demonstrating the weak equidistribution of $\phi(n)$ to moduli $q \le (\log x)^K$ coprime to $6$, obtaining a strong error term that can be expected to be (essentially) of the correct order of magnitude.
\begin{thm}\label{thm:phicoprime6}
Fix $K \ge 1$ and $\epsilon \in (0, 1)$.  We have
\begin{multline*}
\#\{n \le x: \phi(n) \equiv a \pmod q\}\\ = \phiqrec \#\{n \le x: \gcd(\phi(n), q)=1\} + O_{K, \epsilon}\left(\frac x{\phi(q)(\log x)^{1-\alpha(q)(1/3+\epsilon)}}\right),
\end{multline*}
uniformly in moduli $q \le (\log x)^K$ that are coprime to $6$ and in coprime residue classes $a$ mod $q$, where $\alpha(q) \coloneqq \prod_{\ell\mid q } \left(1-\frac1{\ell-1}\right)$. In particular, $\phi(n)$ is weakly equidistributed to moduli $q$ coprime to $6$ growing uniformly up to $(\log x)^K$.  
\end{thm} 
The ``$1/3$" in the exponent of $\log x$ in the error term of the estimate above arises as the maximum absolute value of the averages $\rho_\chi \coloneqq \phiqrec \sum_{\substack{v \bmod q\\(v, q)=1}} ~ \chi(v-1)$ taken over the nontrivial Dirichlet characters $\chi$ mod $q$. The averages $\rho_\chi$ play the roles of the parameter $\varrho$ in the corresponding analogues of Theorem \ref{thm:mainbound} when the role of the multiplicative function $f$ is played by the functions $\chi \circ \phi$. These functions in turn arise by applying the orthogonality of Dirichlet characters modulo $q$ to detect the congruence $\phi(n) \equiv a \pmod q$. 
As shown in section \ref{sec:Thm1.31.4Proof}, the maximum value of $|\rho_\chi|$ is attained by the characters $\chi$ having conductor $5$ (in the case when $q$ is divisible by $5$). This suggests (but does not prove) that the constant ``$1/3$" in the error term of Theorem \ref{thm:phicoprime6} cannot be replaced by a smaller constant in general.

We remark that in the course of our arguments for Theorem \ref{thm:phicoprime6} (as well as Theorem \ref{thm:phidiv3} below), we shall obtain sharp upper bounds for the character sums $\chi(\phi(n))$ for nontrivial Dirichlet characters $\chi$ to our moduli $q$. In the range $q \le (\log x)^K$, these estimates significantly improve upon those given in the work of Balasuriya, Shparlinski and Sutantyo (the case $f(n)=1$ in  \cite[Theorem 1]{BSS09}; see \cite{BS04} and \cite{BS06} for related results on exponential sums involving the Euler totient).

In a forthcoming manuscript \cite{SR23}, the second named author establishes a variant of Theorem \ref{thm:mainbound}, that is useful in studying the distribution of the sum-of-divisors function $\sigma(n) \coloneqq \sum_{d \mid n} d$ among the coprime residue classes to moduli varying uniformly in a wide range.

We also give an analogue of Theorem \ref{thm:phicoprime6} for moduli $q$ that are divisible by $3$. Work of Dence and Pomerance \cite[Theorem 3.1]{DP98} shows that for $i \in \{1, 2\}$, 
$$\#\{n \le x: \phi(n) \equiv i \pmod 3\} \sim c_i \frac{x}{\sqrt{\log x}} \quadxtoinfty,$$
where $c_1 \approx 0.6109\dots$ and $c_2 \approx 0.3284\dots$. 
Our next result shows a similar phenomenon for odd moduli divisible by $3$ that vary uniformly up to a fixed power of $\log x$: we show that the count of $n \le x$ for which $\phi(n)$ lies in a given reduced residue class $a$ mod $q$ is determined by the count of $n \le x$ having $\phi(n) \equiv a \pmod 3$.
\begin{thm}\label{thm:phidiv3}
Fix $K \ge 1$ and $\epsilon \in (0, 1)$.  We have
\begin{multline}\label{eq:phidiv3Est1}
\#\{n \le x: \phi(n) \equiv a \pmod q\} = \frac2{\phi(q)} \#\{n \le x: \gcd(\phi(n), q)=1, \phi(n) \equiv a \pmod 3\}\\ + O_{K, \epsilon}\left(\frac x{\phi(q)(\log x)^{1-\alpha(q)(1/3+\epsilon)}}\right),
\end{multline}
uniformly in coprime residue classes $a$ to moduli $q \le (\log x)^K$ satisfying $\gcd(q, 6)=3$, where $\alpha(q) \coloneqq \prod_{\ell\mid q } \left(1-\frac1{\ell-1}\right)$. As a consequence, 
\begin{equation}\label{eq:phidiv3Est2}
\#\{n \le x: \phi(n) \equiv a \pmod q\} \sim \frac2{\phi(q)} \#\{n \le x: \gcd(\phi(n), q)=1, \phi(n) \equiv a \pmod 3\}
\end{equation}
\text{ as }$x \rightarrow \infty$, in the same range of uniformity in $q$ and $a$. 
\end{thm}
We remark that it is possible to give fairly precise estimates for the main terms in Theorems \ref{thm:phicoprime6} and \ref{thm:phidiv3} by means of Theorems A and B in \cite{scourfield85} (see the remark at the end of section \ref{sec:Thm1.31.4Proof}). See also \cite{scourfield84} and \cite[Proposition 2.1]{PSR23} for similar estimates.

\subsection*{Notation and conventions:} We do not consider the zero function as multiplicative (thus, if $f$ is multiplicative, then $f(1)=1$). We write $P(n)$ or $P^+(n)$ for the largest prime divisor of $n$, with the convention that $P(1)=1$. We set $P_1(n) \coloneqq P(n)$ and define, inductively, $P_k(n) \coloneqq P_{k-1}(n/P(n))$; thus, $P_k(n)$ is the $k$th largest prime factor of $n$ (counted with multiplicity), with $P_k(n)=1$ if $\Omega(n) < k$.  When there is no danger of confusion, we write $(a,b)$ instead of $\gcd(a,b)$. Throughout, the letters $p$ and $\ell$ are to be understood as denoting primes. Implied constants in $\ll$ and $O$-notation may always depend on any parameters declared as ``fixed''; other dependence will be noted explicitly (for example, with subscripts). We write $\log_{k}$ for the $k$th iterate of the natural logarithm.
\section{Uniform bounds on the partial sums of multiplicative functions: Proof of Theorem \ref{thm:mainbound}}
We start by removing the $n \le x$ that are either $z$-smooth or have a repeated prime factor exceeding $y$. Since $|f(n)| \le 1$, the contribution of the former $n$ is bounded in absolute value by $\Psi(x, z)$, while that of the latter $n$ is absolutely bounded by
$$\sum_{p>y} \sum_{\substack{n \le x\\p^2\mid n}} 1 \le x\sum_{p>y} \frac1{p^2} \ll \frac xy,$$
which is also absorbed in the expressions given in the claimed bounds. 

It remains to deal with the sum of $f(n)$ over the surviving $n$, namely those that have $P^+(n)>z$ and no repeated prime factor exceeding $y$; we denote this sum by $\sumxyzf$. Such $n$ can be uniquely written in the form $m P_j \cdots P_1$ for some $j \ge 1$, where $P_1 \coloneqq P^+(n) > z$ and $P^+(m) \le y < P_j < \dots < P_1$. As such $f(n) = f(m) f(P_j) \cdots f(P_1)$ and 
\begin{align*}
\sumxyzf &= \sum_{j \ge 1} \summ f(m) \sum_{\substack{P_1, \dots, P_j\\ P_1>z,~ P_j \cdots P_1 \le x/m\\y < P_j < \dots < P_1}} f(P_j) \cdots f(P_1)\\
&= \sumj \summ f(m) \sum_{\substack{P_2, \dots, P_j\\P_j \cdots P_2 \le x/mz\\y<P_j < \dots < P_2\\}} f(P_j) \cdots f(P_2) \sum_{\max\{P_2, z\}<P_1 \le x/m P_2 \cdots P_j} f(P_1).
\end{align*}
Here $\max\{P_2, z\}$ is to be replaced by ``$z$" in the case $j=1$.

Invoking the hypothesis \eqref{eq:mainhypo} for the innermost sum on $P_1$, the total size of the resulting error term is 
$$\ll Mx \err(y) \sumj \sum_{\substack{m, P_2, \dots, P_j\\mP_2 \cdots P_j \le x/z\\P^+(m) \le y < P_j < \dots < P_2}} \frac1{mP_2 \cdots P_j} \le Mx \err(y) \sum_{n \le x/z} \frac1n \ll Mx\err(y) \log x.$$
We now obtain
\begin{equation}\label{eq:PostP1Removal}
\begin{split}
\sumxyzf = \varrho\sumj \frac1{(j-1)!} \summ f(m) \sum_{z<P_1 \le x/m} ~ 
\sumPTwoj &f(P_2) \cdots f(P_j)\\ &+ O(Mx\err(y) \log x). 
\end{split}   
\end{equation}
Now for $j \ge 2$ and each $i \in \{2, \dots, j\}$, the hypothesis \eqref{eq:mainhypo} yields
\begin{equation}\label{eq:SumfPi}
\begin{split}
\sum_{\substack{y<P_i \le x/mP_1 \cdots P_{i-1}P_{i+1} \cdots P_j \\ P_i<P_1, ~ r \ne i \implies P_r \ne P_i}} f(P_i) &= \sum_{\substack{y<P_i \le x/mP_1 \cdots P_{i-1}P_{i+1} \cdots P_j\\P_i < P_1}} f(P_i) + O(j)\\
&= \varrho\sum_{\substack{y<P_i \le x/mP_1 \cdots P_{i-1}P_{i+1} \cdots P_j\\P_i<P_1, ~ r \ne i \implies P_r \ne P_i}} 1 + O\left(j + \frac{Mx}{mP_1 \cdots P_{i-1}P_{i+1} \cdots P_j}  \err(y)\right).
\end{split}
\end{equation}
In order to estimate the innermost sum in the main term of \eqref{eq:PostP1Removal}, we invoke this estimate successively for all $j \ge 2$ and each $i \in \{2, \dots, j\}$. Indeed 
\begin{align*}
\sumPTwoj f(P_2) \cdots f(P_j) &= \sum_{\substack{P_3, \dots, P_j \in (y, P_1)\\P_3, \dots, P_j \text{ distinct}\\P_3 \cdots P_j \le x/myP_1}} f(P_3) \cdots f(P_j) \sum_{\substack{y<P_2 \le x/mP_1 P_3 \cdots P_j \\ P_2<P_1, ~ r \ne 2 \implies P_2 \ne P_r}} f(P_2)\\
&= \varrho \sum_{\substack{P_2, \dots, P_j \in (y, P_1)\\P_2, \dots, P_j \text{ distinct}\\P_2 \cdots P_j \le x/mP_1}} f(P_3) \cdots f(P_j) + O(\widetilde{\err}),
\end{align*}
where 
$$\widetilde{\err} \coloneqq j \sumPTwojMO 1 + \frac{Mx \err(y)}{mP_1} \sumPTwojMO \frac1{P_2 \cdots P_{j-1}},$$
and we have noted that the error term resulting from the application of \eqref{eq:SumfPi} for $i=2$ is, by relabelling, equal to $\widetilde{\err}$. Likewise, invoking \eqref{eq:SumfPi} for $i = 3, \dots, j$ to successively remove the $f(P_3), \dots, f(P_j)$, we obtain 
$$\sumPTwoj f(P_2) \cdots f(P_j) = \varrho^{j-1} \sumPTwoj 1 + O(j \widetilde{\err}).$$
Plugging this into \eqref{eq:PostP1Removal} for each $j \ge 2$, we obtain an error term of size 
\begin{align*}\allowdisplaybreaks
&\ll \sum_{j \ge 2} \frac1{(j-1)!} \summ \sum_{z<P_1 \le \frac{x}{my^{j-1}}} \Bigg\{j^2 \sumPTwojMO 1 + \frac{Mxj \err(y)}{mP_1} \sumPTwojMO \frac1{P_2 \cdots P_{j-1}}\Bigg\}\\
&\ll (\log x)^2 \sum_{j \ge 2} \sum_{\substack{m, P_1, \dots, P_{j-1}\\mP_1 \cdots P_{j-1} \le x/y\\P_1>z, ~ P^+(m) \le y< P_{j-1} < \dots <P_1}} 1 + Mx (\log x) \err(y) \sum_{\substack{m, P_1, \dots, P_{j-1}\\mP_1 \cdots P_{j-1} \le x/y\\P_1>z, ~ P^+(m) \le y< P_{j-1} < \dots <P_1}} \frac1{m P_1 \cdots P_{j-1}}\\
&\ll (\log x)^2 \sum_{n \le x/y} 1 ~ + ~ M x (\log x) \err(y) \sum_{n \le x/y} \frac1n \ll Mx (\log x)^2 \left(\err(y) + \frac1y\right), 
\end{align*}
leading to the estimate
\begin{equation}\label{eq:PostP1PjRemoval}
\sumxyzf = \varrho\sumj \frac{\varrho^{j-1}}{(j-1)!} \summ f(m) \sum_{\substack{P_1, \dots, P_j\\P_1>z; ~ P_1 \cdots P_j \le x/m\\P_2 \cdots P_j \in (y, P_1) \text{ distinct}}} 1 + O\left(Mx (\log x)^2 \left(\err(y) + \frac1y\right)\right).
\end{equation}
Now the main term above is absolutely bounded by 
\begin{multline*}
|\varrho|\sumj \frac{|\varrho|^{j-1}}{(j-1)!} \summ |f(m)| \sum_{\substack{P_2, \dots, P_j \in (y, x)\\P_2 \cdots P_j \le x/mz}} ~ \sum_{z < P_1 \le x/mP_2 \cdots P_j} 1\\
\ll \frac{|\varrho| x}{\log z}\sumj \frac{|\varrho|^{j-1}}{(j-1)!} \summ \frac{|f(m)|}m \sum_{\substack{P_2, \dots, P_j \in (y, x)}} \frac1{P_2 \cdots P_j}.
\end{multline*}
Since $|f(n)|\le 1$ for all $n$, the sum on $m$ above is no more than 
$$\prod_{p \le y} \left(1+\sum_{v \ge 0} \frac{|f(p^v)|}{p^v}\right) \le \exp\left(\sum_{p \le y} \frac{|f(p)|}p + O\left(\sum_{p \le y}\frac1{p^2}\right)\right) \ll \exp\left(\sum_{p \le y} \frac{|f(p)|}p\right).$$
On the other hand, the sum on $P_2, \dots, P_j$ is 
$$\le \left(\sum_{y<p \le x} \frac1p\right)^{j-1} \le \left(\log\left(\frac{\log x}{\log y}\right) + O(1)\right)^{j-1}.$$
Collecting estimates, we find that the main term in \eqref{eq:PostP1PjRemoval} is 
\begin{multline*}
\ll \frac{|\varrho| x}{\log z} \exp\left(\sum_{p \le y} \frac{|f(p)|}p\right) \sumj \frac1{(j-1)!} \left(|\varrho| \log\left(\frac{\log x}{\log y}\right) + O(1)\right)^{j-1}\\
\ll \frac{|\varrho| x}{\log z} \left(\frac{\log x}{\log y}\right)^{|\varrho|} \exp\left(\sum_{p \le y} \frac{|f(p)|}p\right),
\end{multline*}
completing the proof of the theorem. 
\section{Equidistribution of the Alladi-Erd\H{o}s function: Proof of Theorem \ref{thm:A(n)Effective}}
By the orthogonality of the additive characters mod $q$, we have
\begin{equation}\label{eq:A(n)Orth}
\sum_{\substack{n \le x\\A(n) \equiv a \pmod q}} 1 = \frac xq + \frac1q \sum_{0<r<q} e\left(-\frac{ar}q\right) \sum_{n \le x} e\left(\frac{rA(n)}q\right).
\end{equation}
Let $y \coloneqq \exp((\log x)^{\epsilon/2})$ and $z \coloneqq x^{1/\log_2 x}$. Then for each $r \in \{0, 1, \dots , q-1\}$, the multiplicative function $f(n) \coloneqq e(rA(n)/q)$ satisfies the hypothesis \eqref{eq:mainhypo} with $\varrho \coloneqq \rho_r \coloneqq \phiqrec \sum_{\substack{v \bmod q\\ (v, q)=1}} e\left(\frac{rv}q\right)$. Indeed by the Siegel-Walfisz theorem, we find that
\begin{equation}\label{eq:erApq}
\begin{split}
\sum_{y<p \le Y} e\left(\frac{rA(p)}q\right) &= \sum_{y<p \le Y} e\left(\frac{rp}q\right) = \sum_{\substack{v \bmod q\\ (v, q)=1}} e\left(\frac{rv}q\right) \sum_{\substack{y<p \le Y\\p \equiv v \pmod q}} 1\\
&= \sum_{\substack{v \bmod q\\ (v, q)=1}} e\left(\frac{rv}q\right) \left\{\phiqrec \sum_{\substack{y<p \le Y}} 1 + O\left(Y \exp(-C_0\sqrt{\log Y})\right)\right\}\\ &= \rho_r(\pi(Y) - \pi(y)) + O(\phi(q) Y \exp(-C_0\sqrt{\log y})),
\end{split}
\end{equation}
uniformly for all $q \le (\log x)^K$ and $Y \ge y$. Here $C_0>0$ is an absolute constant and we have noted that $(\log x)^K < y/2$ for all sufficiently large $x$. This verifies hypothesis \eqref{eq:mainhypo} with $\varrho \coloneqq \rho_r$, $M \coloneqq \phi(q)$ and $\err(y) \coloneqq \exp(-C_0\sqrt{\log y})$. From Theorem \ref{thm:mainbound}, 
we obtain 
\begin{equation}\label{eq:e(rA(n))/q}
\sum_{n \le x} e\left(\frac{rA(n)}q\right) \ll \frac{|\rho_r| x}{(\log x)^{1-|\rho_r|-2\epsilon/3}} + \frac x{(\log x)^{(1+o(1)) \log_3 x}},
\end{equation}
where we have noted that $\Psi(x, z) \ll x/(\log x)^{(1+o(1)) \log_3 x}$ by well-known results on smooth numbers (for instance \cite[Theorem 5.13\text{ and }Corollary 5.19, Chapter III.5]{tenenbaum15}).

The sum $\rho_r$ is a (normalized) Ramanujan sum and is nonzero only if $q'_r \coloneqq q/(q, r)$ is squarefree (see \cite{HW08} or \cite{MV07} for some standard properties of Ramanujan sums). Moreover, unless $q$ is even and $r=q/2$, we have $q'_r \ge 3$, so that $|\rho_r| \le \frac1{\phi(q'_r)} \le 1/2$. This shows that 
\begin{equation}\label{eq:Allbutq/2}
\sum_{\substack{0<r<q\\r \ne q/2}} \left|\sum_{n \le x} e\left(\frac{rA(n)}q\right) \right| \ll \Bigg(\sum_{\substack{0<r<q\\q'_r \text{ squarefree }}} \frac1{\phi(q'_r)}\Bigg) \frac x{(\log x)^{1/2-2\epsilon/3}} + \frac x{(\log x)^{\frac12\log_3 x}} \ll \frac x{(\log x)^{1/2-\epsilon}};
\end{equation}
here we have noted that for each squarefree divisor $d$ of $q$, there are exactly $\phi(d)$ many residues $r$ mod $q$ for which $q'_r = q/(q, r) = d$, so that 
$$\sum_{\substack{0<r<q\\q'_r \text{ squarefree }}} \frac1{\phi(q'_r)} \le \sum_{\substack{d\mid q \\d\text{ squarefree}}} 1 = 2^{\omega(q)} \le (\log x)^{\epsilon/6}.$$
It remains to deal with the case $r = q/2$, which arises only for even $q$. But for this value of $r$, a classical estimate of Hall and Tenenbaum for the mean values of multiplicative functions taking values in $[-1, 1]$ (see \cite[Theorem 4.14]{tenenbaum15}) yields
\begin{equation}\label{eq:q/2AddvCharac}
\sum_{n \le x} e\left(\frac{rA(n)}q\right) = \sum_{n \le x} (-1)^{A(n)} \ll \frac x{(\log x)^{3/5}}.
\end{equation}
(In fact, Alladi and Erd\H{o}s \cite{alladi77} show that the left hand side of \eqref{eq:q/2AddvCharac} is $O(x \exp(-c_0 \sqrt{\log x \log_2 x}))$ for some absolute constant $c_0>0$.) Inserting \eqref{eq:Allbutq/2} and \eqref{eq:q/2AddvCharac} into \eqref{eq:A(n)Orth} completes the proof of the theorem.

\section{Distribution of Euler's totient to odd moduli:\\ Proof of Theorems \ref{thm:phicoprime6} and \ref{thm:phidiv3}} \label{sec:Thm1.31.4Proof}
In what follows, we assume $q$ to be any odd positive integer less than or equal to $(\log x)^K$, until stated otherwise. We abbreviate $\alpha(q)$ to $\alpha$; it will be useful to note that for all odd $q$, we have $\alpha \gg 1/\log_2 (3q)$. We first state a rough estimate on the count of $n \le x$ having $\phi(n)$ coprime to $q$, for odd numbers $q \le (\log x)^K$. 
\begin{prop}\label{prop:phiqcoprim}
Fix $K>0$. We have 
\begin{equation*}
\#\{n\le x: (\phi(n),q)=1\} = \frac{x}{(\log{x})^{1-\alpha}} \exp(O((\log\log{(3q)})^{O(1)})), \end{equation*}
uniformly in odd $q \le (\log x)^K$. 
\end{prop}
The above result is a special case of Proposition 2.1 in \cite{PSR23} (cf. work of Scourfield \cite{scourfield84, scourfield85} who gives more precise estimates). In view of Proposition \ref{prop:phiqcoprim}, it remains to show the first assertion of Theorem \ref{thm:phicoprime6} in order to complete the proof of the theorem. To that end, we also state the following lemma, which is a special case of \cite[Lemma 2.4]{PSR23}.   
\begin{lem}\label{lem:primesumphi} For each positive integer $q$ and each real number $x\ge 3q$, \[ \sum_{p \le x} \frac{\1_{(p-1,q)=1}}{p} = \alpha \log_2{x} + O((\log\log{(3q)})^{O(1)}).\]
\end{lem} 
Coming to the first assertion of Theorem \ref{thm:phicoprime6}, set $y \coloneqq \exp((\log x)^{\epsilon/2})$ and $z \coloneqq x^{1/\log_2 x}$ as in the proof of Theorem \ref{thm:A(n)Effective}. We start by removing from the count of $n \le x$ satisfying $\phi(n) \equiv a \pmod q$ those that are either $z$-smooth or have a repeated prime factor exceeding $y$. As observed before, the total contribution of such $n$ is $\ll \Psi(x, z) + x/y \ll x/(\log x)^{(1+o(1))\log_3 x}$, which is negligible in comparison to the error term in the statement of the theorem. 

Among the surviving $n$, we also remove those that have $P_2(n) \le y$: any such $n$ can be written in the form $n = mP$, where $P = P(n)>z$, $m$ is $y$-smooth and $\phi(n) = \phi(m) \phi(P) = \phi(m) (P-1)$. As such, the condition $\phi(n) \equiv a \pmod q$ forces $\phi(m)$ to be coprime to $q$ and, for each choice of $m$, constrains $P \in (z, x/m]$ to at most one coprime residue class modulo $q$. By the Brun-Titchmarsh theorem, there are $\ll x/\phi(q)m\log(z/q) \ll x\log_2 x/\phi(q)m\log x$ many possible choices of $P$ for each choice of $m$. Consequently, the total contribution of the surviving $n \le x$ having $P_2(n) \le y$ is
\begin{equation*}
\begin{split}
\ll \frac{x\log_2 x}{\phi(q)\log x}\sum_{\substack{m: ~ P^+(m) \le y}} \frac{\bbm_{(\phi(m), q)=1}}m &\ll \frac{x\log_2 x}{\phi(q)\log x} \exp\left(\sum_{p \le y} \frac{\bbm_{(p-1,  q)=1}}p\right)\\ 
&\ll \frac{x\log_2 x}{\phi(q)(\log x)^{1-\alpha\epsilon/2}} \exp((\log_2 (3q))^{O(1)}) \ll \frac{x}{\phi(q)(\log x)^{1-2\alpha\epsilon/3}},
\end{split}
\end{equation*}
where we invoked Lemma \ref{lem:primesumphi} in order to estimate the sum $\sum_{p \le y} \bbm_{(p-1, q)=1}/p$ and recalled that $\alpha \gg 1/\log_3 x$ for all odd $q \le (\log x)^K$. Collecting estimates, we have so far shown that
\begin{equation*}
\sum_{\substack{n \le x\\\phi(n) \equiv a \pmod q}} 1 = \sum_{\substack{n \le x\\P(n)>z, ~ P_2(n)>y\\p>y \implies p^2 \nmid n}} \bbm_{\phi(n) \equiv a \pmod q} + O\left(\frac{x}{\phi(q)(\log x)^{1-2\alpha\epsilon/3}}\right).
\end{equation*}
By the orthogonality of the Dirichlet characters mod $q$, we have $\bbm_{\phi(n) \equiv a \pmod q} = \frac{\bbm_{(\phi(n), q)=1}}{\phi(q)} + \phiqrec\sumnontrivchar \chibara \chiphin$, where $\chi_0$ denotes the principal character modulo $q$, and the last sum is over the nonprincipal Dirichlet characters $\chi$ mod $q$. This yields 
\begin{equation}\label{eq:Postfirstremovals}
\begin{split}
\sum_{\substack{n \le x\\\phi(n) \equiv a \pmod q}} 1 = \phiqrec&\sum_{\substack{n \le x\\(\phi(n), q)=1}} 1\\ &+ \phiqrec\sumnontrivchar \chibara\sum_{\substack{n \le x\\P(n)>z, ~ P_2(n)>y\\p>y \implies p^2 \nmid n}} \chiphin + O\left(\frac{x}{\phi(q)(\log x)^{1-2\alpha\epsilon/3}}\right).
\end{split}
\end{equation}
Here we have used the same arguments as before to see that there are $O({x}/{(\log x)^{1-2\alpha\epsilon/3}})$ many $n \le x$ satisfying $(\phi(n), q)=1$ but failing at least one of the following three conditions:
\begin{itemize}
    \item[(i)] $P(n)>z,$
    \item[(ii)] $p>y \implies p^2 \nmid n,$ 
    \item[(iii)] $P_2(n)>y$.  
\end{itemize}
Indeed, any $n$ satisfying conditions (i) and (ii) but failing condition (iii) is of the form $mP$ with $P^+(m) \le y$, $P = P^+(n) \in (z, x/m]$ and $\phi(n) = \phi(m)(P-1)$. As such, $(\phi(m), q)=1$ and the number of $P$ given $m$ is $\ll x/m\log z \ll x\log_2 x/m\log x$, summing which over $m$ yields the claimed bound. 

We now apply the methods from the proof of Theorem \ref{thm:mainbound} in order to estimate the inner sums of $\chi(\phi(n))$ occurring in \eqref{eq:Postfirstremovals}. Any $n$ having $P(n)>z$, $P_2(n)>y$ and no repeated prime factor exceeding $y$ can be uniquely written in the form $mP_j \cdots P_1$ for some $j \ge 2$, where $P_1 = P(n) > z$ and $P(m) \le y < P_j < \dots < P_1$. This shows that 
$$\sum_{\substack{n \le x\\P(n)>z, ~ P_2(n)>y\\p>y \implies p^2 \nmid n}} \chiphin = \sum_{j \ge 2} \sum_{\substack{m \le x\\P(m) \le y}} \chi(\phi(m)) \sum_{\substack{P_1, \dots, P_j\\P_j \cdots P_1 \le x/m\\P_1>z, ~ y<P_j< \dots < P_1}} \chi(P_1-1) \cdots \chi(P_j-1).$$
We proceed as in the proof of Theorem 1 to successively remove $\chi(P_1-1), \cdots, \chi(P_j-1)$, with the input from \eqref{eq:mainhypo} replaced by the estimate for $\sum_{y < p \le Y} \chi(p-1)$ analogous to \eqref{eq:erApq} and coming from the Siegel-Walfisz Theorem. After the dust settles, we are left with 
\begin{equation*}
\begin{split}
\sum_{\substack{n \le x\\P(n)>z, ~ P_2(n)>y\\p>y \implies p^2 \nmid n}} \chiphin = \sum_{j \ge 2} \frac{(\rho_\chi)^j}{(j-1)!} \sum_{\substack{m \le x\\P(m) \le y}} \chi(\phi(m)) &\sum_{\substack{P_1, \dots, P_j\\P_1>z, ~ P_j \cdots P_1 \le x/m\\P_2, \dots , P_j \in (y, P_1) \text{ distinct }}} 1\\ &+ O(x \exp(-K_0 \sqrt{\log y})),
\end{split} 
\end{equation*}
where $\rho_\chi \coloneqq \phiqrec \sum_{v \bmod q} \chi_0(v) \chi(v-1)$ and $K_0>0$ is a constant depending at most on $K$. Hereafter, carrying out the rest of the proof of Theorem \ref{thm:mainbound} yields
\begin{equation}\label{eq:chisign_Oddq1}
\begin{split}
\sum_{\substack{n \le x\\P(n)>z, ~ P_2(n)>y\\p>y \implies p^2 \nmid n}} \chiphin &\ll |\rho_\chi|^2 ~ \frac{x \log_2 x}{\log z} \cdot (\log x)^{|\rho_\chi|} ~ \exp\left(\sum_{p \le y} \frac{\bbm_{(p-1, q)=1}}p \right) + x \exp(-K_0 \sqrt{\log y})\\
&\ll |\rho_\chi|^2 \frac x{(\log x)^{1-|\rho_\chi|-2\alpha\epsilon/3}} + x \exp(-K_0 \sqrt{\log y}),
\end{split}
\end{equation}
where we have recalled Lemma \ref{lem:primesumphi} while passing to the second line above. 

To use the last bound, we need to better understand the $|\rho_\chi|$. Let $\condchi$ denote the conductor of $\chi \bmod q$, so that $\condchi>1$ and $\condchi \mid  q$. We can write $\chi$ uniquely in the form $\prod_{\ell^e \parallel q} \chi_\ell$, where each $\chi_\ell$ is a character mod $\ell^e$ and $\chi_\ell$ is nontrivial precisely when $\ell\mid \condchi$. Now $\phi(q)\rho_\chi = \prod_{\ell^e \parallel q} S_{\chi, \ell}$, where for each prime power $\ell^e \parallel q$, 
\begin{equation}\label{eq:Schiell}
S_{\chi, \ell} \coloneqq \sum_{v \bmod \ell^e} \chi_{0, \ell}(v) \chi_\ell(v-1) = \sum_{\substack{v \bmod{\ell^e}\\(v, \ell)=1}} \chi_\ell(v-1) = \sum_{u \bmod{\ell^e}} \chi_\ell(u) - \sum_{\substack{u \bmod{\ell^e}\\u \equiv -1 \pmod \ell}} \chi_\ell(u).
\end{equation}
Here $\chi_{0, \ell}$ denotes the trivial character mod $\ell^e$, and we have noted that as $v$ runs over the coprime residues mod $\ell^e$, the expression $v-1$ runs over all the residues mod $\ell^e$ except for those congruent to $-1$ mod $\ell$. The first sum in the rightmost expression in \eqref{eq:Schiell} is $\bbm_{\chi_\ell = \chi_{0, \ell}} ~ \phi(\ell^e)$. We claim that the second sum is zero unless $\cond(\chi_\ell) \mid \ell$, in which case it is $\chi_\ell(-1) \ell^{e-1}$. Indeed, the second sum is equal to $\chi_\ell(-1) \sum_{\substack{u \bmod \ell^e\\u \equiv 1 \pmod\ell}} \chi_\ell(u)$, and this latter sum is invariant upon multiplication by an element lying in the subgroup of residues that are $1$ mod $\ell$, and hence is non-vanishing precisely when $\chi_\ell$ restricts to the trivial character on this subgroup. 
This shows that if the second sum is nonzero, then $\cond(\chi_\ell) \mid \ell$. Conversely, if $\cond(\chi_\ell) \mid \ell$, then it is clear that the second sum is equal to $\chi_\ell(-1) \sum_{\substack{u \bmod \ell^e\\u \equiv -1 \pmod \ell}} 1 = \chi_\ell(-1) \ell^{e-1}$, establishing our claim.  

Altogether, we obtain
$$S_{\chi, \ell} = \bbm_{\chi_\ell = \chi_{0, \ell}} ~ \phi(\ell^e) - \bbm_{\cond(\chi_\ell)\mid \ell} ~ \chi_\ell(-1) \ell^{e-1} = \bbm_{\cond(\chi_\ell)\mid \ell} ~ \ell^{e-1}\left(\bbm_{\ell \nmid \condchi} (\ell-1) - \chi_\ell(-1)\right)$$
for each prime power $\ell^e \parallel q$. Multiplying this relation over all these prime powers, we obtain
\begin{equation}\label{eq:rhochi}
\rho_\chi = \prod_{\ell^e \parallel q} \frac{S_{\chi, \ell}}{\phi(\ell^e)} = \bbm_{\condchi \text{ squarefree }} \prod_{\ell^e \parallel q} \left(\bbm_{\ell \nmid \condchi} - \frac{\chi_\ell(-1)}{\ell-1}\right) = \bbm_{\condchi \text{ squarefree }} \frac{(-1)^{\omega(\condchi)} \chi(-1) \alpha}{\prod_{\ell\mid \condchi} (\ell-2)}.
\end{equation}
If $3\mid q$, let $\psi$ denote the unique character mod $q$ induced by the nontrivial character mod $3$. Then for any nontrivial character $\chi \ne \psi$ mod $q$ for which $\rho_\chi \ne 0$, its conductor $\condchi$ has a prime divisor at least $5$, so that $|\rho_\chi| \le \alpha/3$ by \eqref{eq:rhochi}. As such, \eqref{eq:chisign_Oddq1} yields for all such $\chi$,
\begin{equation}\label{eq:chisign_Oddq2_NonExcp}
\sum_{\substack{n \le x\\P(n)>z, ~ P_2(n)>y\\p>y \implies p^2 \nmid n}} \chiphin \ll |\rho_\chi|^2 \frac x{(\log x)^{1-(1/3+2\epsilon/3)\alpha}} + x \exp(-K_0 \sqrt{\log y}). 
\end{equation}
But from \eqref{eq:rhochi} and the fact that there are exactly $\prod_{\ell\mid d} (\ell-2)$ primitive characters modulo any squarefree integer $d$, we find that
\begin{equation*}
\begin{split}
\sum_{\chi \bmod q} |\rho_\chi|^2 \le \alpha^2\sum_{\substack{d\mid q\\d\text{ squarefree}}} \frac1{\prod_{\ell\mid d} (\ell-2)^2} \sum_{\substack{\chi \bmod q\\\condchi = d}} 1 \le \alpha^2\sum_{\substack{d\mid q\\d\text{ squarefree}}} \frac1{\prod_{\ell\mid d} (\ell-2)} \le \alpha^2 \prod_{\ell\mid q} \frac{\ell-1}{\ell-2} = \alpha.
\end{split}
\end{equation*}
Summing the bound \eqref{eq:chisign_Oddq2_NonExcp} over all nontrivial characters $\chi \ne \psi$ mod $q$, and plugging the resulting bound into \eqref{eq:Postfirstremovals}, we thus obtain 
\begin{equation}\label{eq:ReducnToExcpChar}
\begin{split}
\sum_{\substack{n \le x\\\phi(n) \equiv a \pmod q}} 1 &= \phiqrec\sum_{\substack{n \le x\\(\phi(n), q)=1}} 1 + \frac{\bbm_{3\mid q} \overline{\psi}(a)}{\phi(q)}\sum_{\substack{n \le x\\P(n)>z, ~ P_2(n)>y\\p>y \implies p^2 \nmid n}} \psi(\phi(n)) + O\left(\frac{x}{\phi(q)(\log x)^{1-\alpha(1/3+\epsilon)}}\right)\\
&= \phiqrec\sum_{\substack{n \le x\\(\phi(n), q)=1}} 1 + \frac{\bbm_{3\mid q} \overline{\psi}(a)}{\phi(q)}\sum_{n \le x} \psi(\phi(n)) + O\left(\frac{x}{\phi(q)(\log x)^{1-\alpha(1/3+\epsilon)}}\right).
\end{split}
\end{equation}
In passing to the second line above, we have recalled our previous bound on the count of $n \le x$ having $\gcd(\phi(n), q)=1$ but failing one of conditions (i)--(iii) in the paragraph following \eqref{eq:Postfirstremovals}.

The last equality in \eqref{eq:ReducnToExcpChar} already establishes the first assertion of Theorem \ref{thm:phicoprime6} for moduli $q \le (\log x)^K$ coprime to $6$, thus completing the proof of that theorem. 

Coming to the proof of Theorem \ref{thm:phidiv3}, we now consider $q \le (\log x)^K$ satisfying $\gcd(6, q)=3$, so that $q$ is odd and divisible by $3$. By definition of $\psi$, we have $$\psi(\phi(n)) = 
\begin{cases}
\psi(a), & \text{ if } \gcd(\phi(n), q)=1 \text{ and } \phi(n) \equiv a \pmod 3,\\
\psi(-a) = -\psi(a), & \text{ if } \gcd(\phi(n), q)=1 \text{ and } \phi(n) \equiv -a \pmod 3,\\
0, & \text{ if } \gcd(\phi(n), q) \ne 1.
\end{cases}$$
Inserting this into \eqref{eq:ReducnToExcpChar} and recalling that $\psi(a) \overline\psi(a) = 1$, we obtain the first estimate claimed in Theorem \ref{thm:phidiv3}. To obtain the final asymptotic formula in the theorem, it suffices to show that
\begin{equation}\label{eq:Reducntox4}
3\sum_{\substack{n \le x: ~ (\phi(n), q)=1\\ \phi(n) \equiv a \pmod 3}} 1 \ge \sum_{\substack{n \le x/4\\ (\phi(n), q)=1}} 1.
\end{equation}
Indeed once we show this, an application of Proposition \ref{prop:phiqcoprim} will reveal that the right hand side of the above inequality grows strictly faster than the $O$-term in \eqref{eq:phidiv3Est1}, which will lead to \eqref{eq:phidiv3Est2}. 

Finally to prove \eqref{eq:Reducntox4}, we split the right hand side of the the inequality as
$$\sum_{\substack{n \le x/4\\ (\phi(n), q)=1\\ \phi(n) \equiv a \pmod 3}} 1 ~ ~ + \sum_{\substack{n \le x/4: ~ 2 \nmid n\\ (\phi(n), q)=1\\ \phi(n) \equiv -a \pmod 3}} 1 ~ ~ + \sum_{\substack{n \le x/4: ~ 2\mid n\\ (\phi(n), q)=1\\ \phi(n) \equiv -a \pmod 3}} 1.$$
We denote the sum on the left hand side of \eqref{eq:Reducntox4} by $S_0$; in other words, $S_0$ is one-third the left hand side of \eqref{eq:Reducntox4}. We claim that each of the three sums in the above display is no more than $S_0$. Indeed, any $n$ counted in the first of the three sums is automatically counted in $S_0$. For any $n$ counted in the second sum above, $\phi(4n) = 2\phi(n)$ is coprime to $q$ (since $q$ is odd) and $\phi(4n) \equiv -\phi(n) \equiv a \pmod 3$, so that $4n \le x$ is counted in $S_0$. Finally, for any $n$ counted in the third sum above, we have $2n \le x/2$, and $\phi(2n) = 2 \phi(n) \equiv -\phi(n) \equiv a \pmod 3$ and $\gcd(\phi(2n), q) = \gcd(2 \phi(n), q)=1$, so that $2n \le x$ is counted in $S_0$. This establishes the bound \eqref{eq:Reducntox4}, completing the proof of Theorem \ref{thm:phidiv3}. 

\begin{rmk}
It is worth pointing out that the count of $n \le x$ satisfying $\gcd(\phi(n), q)=1$ and $\phi(n) \equiv a \pmod 3$ (which appears in the main term on the right hand side of \eqref{eq:phidiv3Est1}) is of the same order of magnitude as the count of $n \le x$ satisfying $\gcd(\phi(n), q)=1$. In other words, these two quantities are bounded by constant multiples of one another (said constants being independent of $q$). This follows by combining the lower bound \eqref{eq:Reducntox4} with Theorems A and B in \cite{scourfield85}. 
\end{rmk}

\providecommand{\bysame}{\leavevmode\hbox to3em{\hrulefill}\thinspace}
\providecommand{\MR}{\relax\ifhmode\unskip\space\fi MR }
\providecommand{\MRhref}[2]{%
  \href{http://www.ams.org/mathscinet-getitem?mr=#1}{#2}
}
\providecommand{\href}[2]{#2}

\end{document}